\documentclass[12pt]{amsart}
\usepackage[dvips]{epsfig}
\usepackage{graphics}
\usepackage{latexsym}
\usepackage{verbatim}
\usepackage{amsmath}
\usepackage{amsthm}
\usepackage{amssymb}
\usepackage{hyperref}
\usepackage{color}
\newcommand{\ens}[1]{\mathbb{#1}}

\newcommand{\N}{\mathbb{N}}

\newcommand{\R}{\mathbb{R}}

\def\cal{\mathcal}

\def\derpar#1#2{\frac{\partial#1}{\partial#2}}
\def\var{\varepsilon}
\def\pa{\partial}

\newcommand{\beqn}{\begin{equation}}
\newcommand{\eeqn}{\end{equation}}
\newcommand{\bear}{\begin{eqnarray}}
\newcommand{\eear}{\end{eqnarray}}
\newcommand{\bean}{\begin{eqnarray*}}
\newcommand{\eean}{\end{eqnarray*}}

\newcommand{\eps}{\varepsilon}

\def\signcm{\bigskip\bigskip\hspace{80mm}
\vbox{{\sc C. Mouhot\par\vspace{3mm}
University of Cambridge, DPMMS \par
Wilberforce road, CB3 0WA, UK \par
\par\vspace{3mm}
e-mail:} C.Mouhot@dpmms.cam.ac.uk}}

\def\signld{\bigskip\bigskip\hspace{80mm}
\vbox{{\sc L. Desvillettes \par\vspace{3mm}
CMLA, ENS Cachan, CNRS, PRES UniverSud \& IUF \par
61, Avenue du Pr\'esident Wilson\par
94235 Cachan Cedex,
FRANCE \par\vspace{3mm}
e-mail:} desville@cmla.ens-cachan.fr}}

\def\signcv{\bigskip\bigskip\hspace{80mm}
\vbox{{\sc C. Villani \par\vspace{3mm}
Institut Henri Poincar\'e \& Universit\'e Claude Bernard Lyon 1\par 
11 rue Pierre et Marie Curie 
75230 Paris Cedex 05, FRANCE  \par\vspace{3mm}
e-mail:} villani@ihp.jussieu.fr }}

\begin{document}

\title[Cercignani's conjecture for the Boltzmann equation] {Celebrating
  Cercignani's conjecture for the Boltzmann equation}

\author{Laurent Desvillettes, Cl\'ement Mouhot \& C\'edric Villani}
\date{\today}

\hyphenation{bounda-ry rea-so-na-ble be-ha-vior pro-per-ties
cha-rac-te-ris-tic}

\begin{abstract}
  Cercignani's conjecture assumes a linear inequality between the
entropy and entropy production functionals for Boltzmann's nonlinear integral operator 
in rarefied gas dynamics. Related to the field of logarithmic Sobolev inequalities
and spectral gap inequalities,
this issue has been at the core of the renewal of the mathematical theory of convergence
to thermodynamical equilibrium for rarefied gases over the past decade.
In this review paper, we survey the various
  positive and negative results which were obtained since the
  conjecture was proposed in the 1980s.
\end{abstract}

\maketitle
{\it This paper is dedicated to the memory of the late Carlo Cercignani, powerful mind and great scientist,
one of the founders of the modern theory of the Boltzmann equation.}

\bigskip



\tableofcontents

\section{The Boltzmann equation and its variants}
\setcounter{equation}{0}

\subsection{The Boltzmann equation}
\label{sec:boltzmann-equation}

After Maxwell \cite{maxwell} wrote down the basic equation for the kinetic theory
of gases, Boltzmann \cite{boltzmann} made such remarkable progress on this equation
that its name remained attached to it. The {\bf Boltzmann equation} describes
the behavior of a rarefied gas when the only interactions taken into
account are binary collisions. This evolution equation reads
\begin{equation*}\label{eq:baseinhom}
 \derpar{f}{t} + v \cdot \nabla_x f = Q(f,f), 
 \qquad  x \in \Omega, \quad v \in \R^d, \quad t \geq 0,
 \end{equation*}
 where $\Omega \subset \R^d$ is the spatial domain ($d\geq 2$) and $f$ is the
 time-dependent particle distribution function in phase space.  
In the case when this distribution function is assumed to be independent of
 the position, equation \eqref{eq:baseinhom} reduces to the 
{\bf spatially homogeneous Boltzmann equation}
 \begin{equation}\label{eq:base}
 \derpar{f}{t}(t,v)  = Q(f,f)(t,v), \quad  v \in \R^d, \quad t \geq 0,
 \end{equation}
 where $Q$ is the quadratic {\bf Boltzmann collision operator},
 defined by the bilinear form
 \begin{equation*}\label{eq:collop}
   Q(g,f) = \int _{\R^d \times \ens{S}^{d-1}} 
   B\bigl(|v-v_*|, \cos \theta\bigr)(g'_* f' - g_* f) \, dv_* \, d\sigma.
 \end{equation*}
We have used the shorthands $f'=f(v')$, $g_*=g(v_*)$ and
$g'_*=g(v'_*)$, where
\begin{equation*}\label{eq:rel:vit}
  v' = \frac{v+v_*}2 + \frac{|v-v_*|}2 \, \sigma  \quad \mbox{ and } \quad 
  v'_* = \frac{v+v_*}2 - \frac{|v-v_*|}2 \, \sigma 
\end{equation*}
stand for the pre-collisional velocities of particles which after collision
have velocities $v$ and $v_*$. Moreover it is usual to introduce the {\bf
  deviation angle} $\theta\in [0,\pi]$ between $v'-v'_*$ and $v-v_*$. The
function $B$ is called the {\bf Boltzmann collision kernel} and it is
determined by physics (it is related to the physical cross-section
$\Sigma(v-v_*,\sigma)$ by the formula $B=|v-v_*| \, \Sigma$). On physical
grounds (in particular Galilean invariance), it is assumed that $B \geq 0$
and $B$ is a function of $|v-v_*|$ and $\cos\theta$ only.

\subsection{The collision kernel}
\label{sec:collision-kernel}

In the theory of Maxwell and Boltzmann, the interaction between particles
is reflected in the formula for the collision kernel $B$.
It may be short-range or long-range. The most important case of
short-range interaction is the {\bf hard spheres} model, where
particles are spheres interacting by contact. In that case,
$B=|v-v_*|$ in dimension $d=3$ (constant cross-section).

Typical models of long-range interactions are given by inverse power-law 
forces \cite{cercignani}. In dimension $d=3$,
if the intermolecular force scales like $r^{-s}$ with $s > 2$, then
$$ B(|v-v_*|, \cos \theta) =  |v-v_*|^\gamma \, b(\theta), \quad \theta \in [0,\pi], $$
where $b$ is smooth except near $\theta =0$,
$$ b(\theta) \sim_{\theta =0} \mbox{cst} \, \theta^{-2 - \nu}, $$
and
$$ \gamma = \frac{s-5}{s-1}, \qquad \nu = \frac{2}{s-1}.$$

It will be interesting to also consider more general
kernels $B$ which do not necessarily come from microscopic interactions;
a remarkable case is $\gamma=2$.

In the important case of Coulomb interaction, the force scales like the inverse of
the square of the distance between particles; then the Boltzmann operator
does not make sense anymore \cite[Annexe~1, Appendice]{villani-hdr}.
This divergence led Landau \cite{landau-1936} to introduce in 1936 a ``diffusive''
version of the Boltzmann collision operator, the {\bf Landau collision operator}, 
defined by the bilinear form
 \begin{multline}\label{eq:colloplandau}
   Q(g,f) = \nabla_v \cdot \Bigg( \int _{\R^N}  \Phi(|v-v_*|)\\
   \bigg\{ |v-v_*|^2\,\mbox{Id} - (v-v_*)\otimes (v-v_*) \bigg\} \, \bigg(
   (\nabla g) f_*  - (\nabla g)_* f
   \bigg)\, dv_*\Bigg)
\end{multline}
leading to the polar form
 \begin{multline}\label{eq:colloplandaupolar} 
   Q(f,f) = \nabla_v \cdot \Bigg( \int _{\R^N}  \Phi(|v-v_*|)\\
     \bigg\{ |v-v_*|^2\,\mbox{Id} - (v-v_*)\otimes (v-v_*) \bigg\}
     \, \bigg( \frac{\nabla f}f - \frac{\nabla f_*}{f_*} \bigg)\, f\,f_*\,
     dv_*\Bigg).
\end{multline}
Here $\Phi(|v-v_*|)= |v-v_*|^{-3}$ and the dimension is $d=3$.
As in \cite{MR1737547,MR1737548} we can consider any dimension $d\geq 2$
and more general functions $\Phi$, say $\Phi(|v-v_*|)=|v-v_*|^\gamma$, $\gamma\geq -d$.

A mathematical way to derive the operator \eqref{eq:colloplandau} is to apply
the {\bf grazing collision limit} to the Boltzmann collision operator with
collision kernel $B=\Phi(|v-v_*|)\,b(\theta)$, that is, to concentrate
on deviation angles $\theta\simeq 0$
\cite{MR1055522,MR1167768,MR1165528,MR1650006,MR2037247}.
The case $\gamma=-3$, $d=3$ considered by Landau in dimension $d=3$ will be called
``Landau-Coulomb'' for the sake of classification.

\subsection{Conserved quantities and entropy structure}

Boltzmann's and Landau's collision operators have the fundamental
properties of conserving mass, momentum and energy
\begin{equation*}
  \int_{\R^d}Q(f,f)(v) \, \phi(v)\,dv = 0, \quad
  \phi(v)=1,v,|v|^2/2 
\end{equation*}
and satisfying (the first part of) Boltzmann's $H$ theorem, which can be
formally written as
\begin{equation*} 
  {\cal D}(f):= - \frac{d}{dt} \int_{\R^d} f \log f \, dv = - \int_{\R^d}
  Q(f,f)\log f \, dv \geq 0.
\end{equation*}
Boltzmann's so-called ``${\cal H}$ functional''
\[
{\cal H}(f) = \int f \, \log f\, dv
\]
is the opposite of the entropy of the gas.

The second part of Boltzmann's $H$ theorem states that under appropriate boundary
conditions, any equilibrium distribution function 
(that is, such that $v\cdot\nabla_x f = Q(f,f)$) satisfies
${\cal D}(f)=0$, or equivalently $Q(f,f)=0$, and takes the form of a Maxwellian distribution
\begin{equation*}
  M(\rho,u,T)(v)=\frac{\rho}{(2\pi T)^{d/2}}
  \exp \left( - \frac{\vert v - u \vert^2} {2T} \right). 
\end{equation*}
The parameters $\rho= 0$, $u\in\R^d$ and $T\geq 0$ are interpreted as respectively
the density, mean velocity and temperature of the gas:
\begin{equation*}
  \rho = \int_{\R^d}f(v) \, dv, \quad u =
  \frac{1}{\rho}\int_{\R^d}v \, f(v) \, dv, \quad T = {1\over{d \rho}}
  \int_{\R^d}\vert v - u \vert^2 \, f(v) \, dv.
\end{equation*}
 As a result of the process of
entropy production pushing towards local equilibrium combined with the
constraints of conservation laws, solutions of the Boltzmann equation
are expected to converge to a unique Maxwellian equilibrium 
(This is the Krasovskii--Lasalle principle in the context of the Boltzmann equation).

Assuming that $\rho$ and $T$ are positive, we may rescale $f:=f(v)$
into $v \mapsto a\, f(\lambda\, (v - b))$, in such a way that the new density,
velocity and temperature are $\rho=1$, $u=0$ and $T=1$.
So we set $M(v)= (2\pi)^{-d/2}\,e^{-|v|^2/2}$ as the Maxwellian equilibrium in the sequel.

\subsection{The linearized collision operators}
\label{sec:line-boltzm-equat}

Let us first consider the Boltzmann collision operator. We introduce the
fluctuation around the Maxwellian equilibrium $M$ computed above: 
$$
f = M + Mh, \quad v \mapsto h(v) \in L^2(M)
$$
where $L^2(M)$ denotes the Lebesgue space $L^2$ on $\R^d$ with reference
measure $M(v) \, dv$. Then the linearized collision operator writes
$$
Lh = M^{-1} \left[ Q(Mh,M) + Q(M,Mh) \right],
$$
or
\begin{equation}\label{eq:LBE}
  Lh = \frac{1}{4} \int_{\R^d \times \R^d \times \ens{S}^{d-1}} 
  \Big[ h' + h'_* - h - h_* \Big] \, B \, M_* \, dv_* \, d\sigma
\end{equation}
for the same collision kernel $B$ as in $Q$.

It is easy to check that $L$ is symmetric in the Hilbert space 
$L^2(M)$ and that it is non-positive in this space (this is the linearized
form of the $H$ theorem). The dissipation of squared $L^2$ norm, that is
the opposite of the Dirichlet form associated with $L$, is
$$
D(h) = - \langle h , Lh \rangle_{L^2(M)} = \frac{1}{4} \int_{\R^d \times
  \R^d \times \ens{S}^{d-1}} \Big| h' + h'_* - h - h_* \Big|^2 \, B \, M \,
M_* \, dv\, dv_* \, d\sigma \ge 0.
$$
It is straightforward from this formula that the null space of $L$ has dimension
$d+2$, and is spanned by the so-called {\bf collisional invariants}
$1,v_1, \dots, v_d, |v|^2$.

In the case of the Landau operator, a similar line of arguments
leads to 
\begin{equation*}
  Lh(v) = M^{-1} \, \nabla _v \cdot \Bigg( \int_{\R^d} {\bf
    a}(v-v_*) \, \Big[ \left( \nabla h \right) - \left( \nabla h \right)_*
  \Big] \, M \, M_* \, dv_* \Bigg),
\end{equation*}
with
\[ {\bf a} (z) = |z|^2 \, \Phi(z) \, \Pi_{z^\bot}, \qquad \left(
  \Pi_{z^\bot} \right) _{i,j} = \delta_{i,j} - \frac{z_i z_j}{|z|^2},
\]
and the corresponding negative Dirichlet form is
\begin{multline}\label{eq:landau-dissip}
D(h) = -\langle h , Lh \rangle_{L^2(M)} \\
= \frac{1}{2} \int_{\R^d \times \R^d} 
\left|v-v_*\right|^2 \, \Phi\left(v-v_*\right) \, 
\Big| \Pi_{(v-v_*)^\bot} \, \Big( \left( \nabla h \right) - \left( \nabla h
    \right)_* \Big) \Big|^2 \, M \, M_* \, dv\, dv_*  \ge 0.
\end{multline}

\subsection{Comparison with usual differential operators and
  classification}
\label{sec:comp-with-usual}

The Boltzmann and Landau collision operators are {\em a priori} extremely
intricate, partly due to their integral or integro-differential nature
(and partly of course due to their nonlinear nature!). Therefore it is
useful, in order to grab an intuition of these operators, to draw a
parallel with usual differential operators  which are more 
familiar.
\par
For the Boltzmann collision operators, say with collision kernel of the
form $B = \Phi(|v-v_*|) \, b(|\theta|)$, the most important two
``parameters'' interplaying and determining its behavior are (1) the growth
or decay of $\Phi$, and (2) the singularity of $b$ at grazing collisions
$\theta \sim 0$. To be more precise, let us consider the model case
$\Phi(z)=z^\gamma$, $\gamma \in (-d,+\infty)$ and $b(|\theta|) \sim
\theta^{-(d-1)-\nu}$, $\nu \in (-\infty,2)$ as $\theta \sim 0$, for the
Boltzmann collision operator. Then the order of singularity (2) plays the
role of the order (highest number of derivatives) in a differential
operator. For instance $\nu <0$ in the model means a zero order operator,
whereas $\nu \in (0,2)$ means a fractional differential operator with order
$\nu$. And the growth or decay of $\Phi$ (1) plays the role of the growth
or decay of the coefficients in a differential operator. Therefore
$\gamma=0$ (the so-called Maxwell or pseudo-Maxwell molecules cases) would
correspond to a constant coefficients differential operator, and $\gamma
=1$ (similar to hard spheres) would correspond to unbounded polynomially
growing coefficients.

From this comparison it becomes natural to consider the Landau
collision operator with $\Phi(z)=z^\gamma$ formally as the limit case
$\nu=2$ in the above classification. This unified picture of this
family of integro-differential operators is summarized in figure
\ref{fig:classification} below.
\begin{figure}[here]
\includegraphics{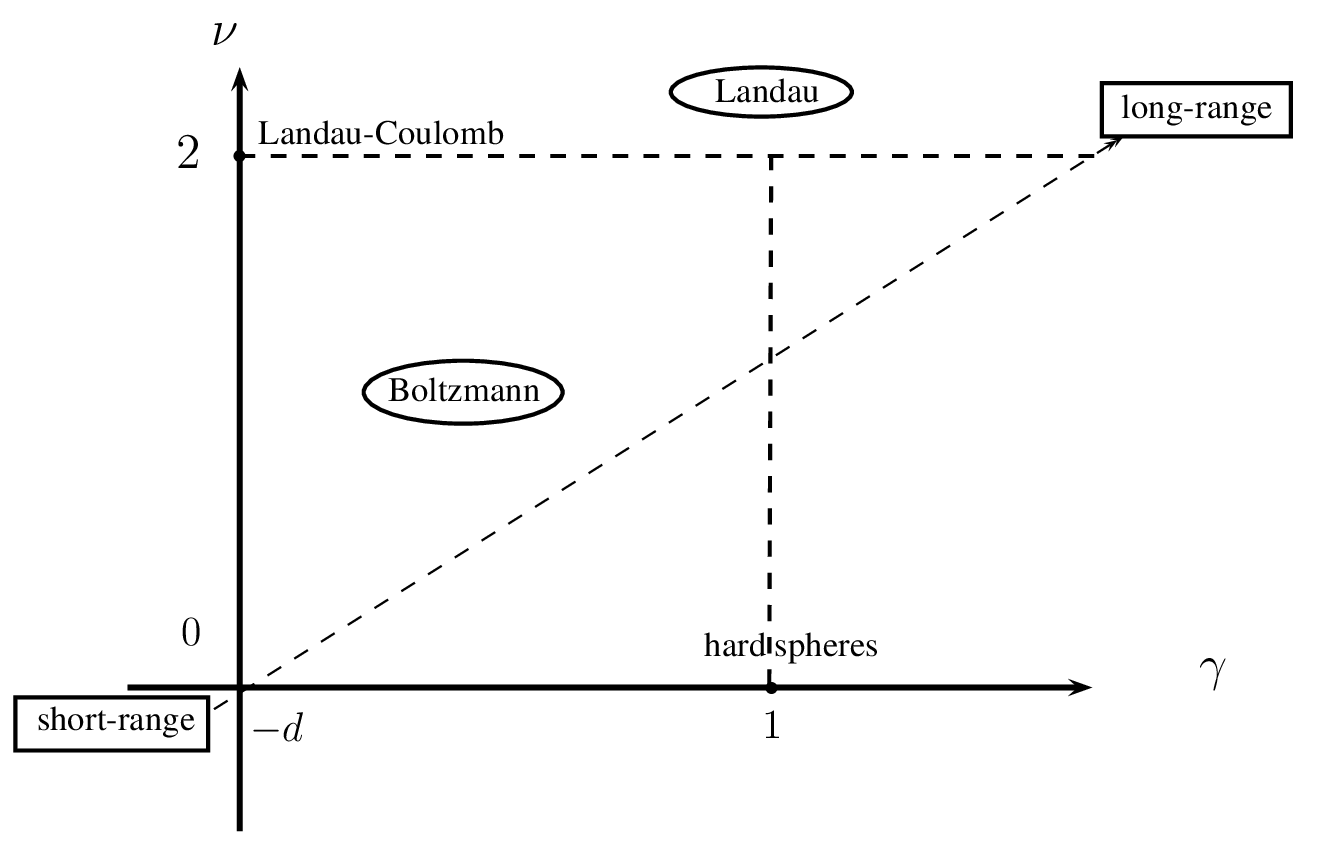}
\caption{Classification of the Boltzmann and Landau operators}
\label{fig:classification}
\end{figure}
\section{Cercignani's conjecture} 
\setcounter{equation}{0}

\subsection{Constructive quantitative estimates for the large time
  behavior}

The relaxation to equilibrium has been studied since the works of Boltzmann and
it is at the core of kinetic theory. The motivation is to provide an
analytic basis for the second principle of thermodynamics for a statistical
physics model of a gas out of equilibrium.  Indeed, Boltzmann's famous
$H$ theorem gives an analytic meaning to the entropy production process and
identifies possible equilibrium states. In this context, proving
convergence towards equilibrium is a fundamental step to justify Boltzmann
model, but cannot be fully satisfactory as long as it remains based on
non-constructive arguments. Indeed, as suggested implicitly by Boltzmann
when answering critics of his theory based on Poincar\'e's recurrence
Theorem, the validity of the Boltzmann equation breaks for very large time
(see~\cite[Chapter~1, Section~2.5]{Vi:hb} for a discussion). It is
therefore crucial to obtain constructive quantitative informations on the
time scale of the convergence, in order to show that this time scale is
much smaller than the time scale of validity of the model. Cercignani's
conjecture is an attempt to provide such constructive quantitative
estimates. In the words of Cercignani: ``{\it The present contribution is
  intended as a step toward the solution of the first main problem of
  kinetic theory, as defined by Truesdell and Muncaster, i.e. ``to discover
  and specify the circumstances that give rise to solutions which persist
  forever''.}'' It is inspired by the entropy - entropy production method,
that we now briefly describe.

\subsection{The entropy - entropy production method}
\label{sec:entropy--entropy}

This method was first used in kinetic theory for the {\bf Fokker-Planck
  equation} (Cf. \cite{BK,TO})
\[ 
\pa_t f = \nabla_v \cdot ( \nabla f + v\, f), \quad v \in \R^d, \quad
\int_{\R^d} f(w) \,dw =1.
\] 
In that case, the equilibrium $M$ is given by the formula
\[ 
M(v) = (2\pi)^{-d/2} \, 
e^{- |v|^2/2}
\] 
and the entropy production is 
\[ 
{\cal D}_{FP}(f) = \int_{\R^d} f(v)\, \bigg| \nabla \log \frac{f}{M} \bigg|^2\,
dv.
\] 
The exponential convergence is then obtained thanks to the {\bf
  logarithmic Sobolev inequality} (cf. \cite{MR0372613}), which
exactly means in this setting
\[ 
{\cal D}_{FP}(f) \ge 2\, \bigg[ {\cal H}(f)  - {\cal H}(M) \bigg].
\]

Consider the more general case of an equation for which a Lyapunov
functional $\mathcal{H}_*$ exists, that is
$$ {\cal D}_{*}(f(t)):= - \frac{d}{dt} {\cal H}_{*}(f(t)) \ge 0 $$ 
and assume that the entropy $-{\cal H}_{*}$ is maximal for a unique
function $M_*$ (among the functions belonging to a space depending on the
conserved quantities in the equation). As seen in the previous section,
this structure is provided by the $H$ theorem for Boltzmann and Landau
equations. The {\bf entropy - entropy production method} consists in
looking for (functional) estimates like
$$
{\cal D}_*(f) \ge \Theta\bigl( {\cal H}_*(f)  - {\cal H}_*(M_*) \bigr), 
$$
where $\Theta : \R_+ \to \R_+$ is a function such that
$$\Theta(x) = 0 \qquad \iff \qquad x=0. $$ 
The stronger $\Theta$ increases near $0$ the better the rate of relaxation
to equilibrium, since the differential inequality
$$  
\frac{d}{dt} \bigl({\cal H}_*(f) - {\cal H}_*(M_*) \bigr) \le -\, 
\Theta\bigl( {\cal H}_*(f) - {\cal H}_*(M_*) \bigr) $$ 
leads to
$$ {\cal H}_*(f(t)) - {\cal H}_*(M_*) \le R(t), $$
where $R$ is the reciprocal of a primitive of {$ - 1/\Theta$}.  Then, if
the relative entropy ${\cal H}_*(f) - {\cal H}_*(M_*)$ is coercive in the
sense that it controls from below some distance or some norm (denoted by
$\|\,\,\|_*$) between $f$ and its associated equilibrium distribution $M_*$
(for the Boltzmann entropy this is precisely provided by the so-called
Czizsar-Kullback-Pinsker inequality, see
\cite{CZ1, CZ2}), we obtain
$$ 
\| f(t) - M_* \|_* \le S(t), 
$$
where (generically)
$$ S(t) =  C \, R(t)^{1/2} . $$

In the particular case $\Phi(x) = C \, x$ (like in the case of the 
Fokker-Planck equation), one gets 
\[
R(t) \le C \, e^{- C' \, t} , 
\] {\em i.e.}, {\bf exponential convergence towards equilibrium} (with
explicit rate). In the slightly worse case $\Phi(x) = C_{\var}\,
x^{1+\var}$ for some (or all) $\var>0$ we can deduce
\[
R(t) \le C_{\var} ' \, t^{-1/\var},
\] 
and we thus get {\bf algebraic convergence towards equilibrium} (with
explicit rate). When $\var>0$ can be taken as small as one wishes, we speak
of {\bf almost exponential convergence}.

\subsection{Cercignani's conjecture}
\label{sec:cerc-conj}

The original Cercignani's conjecture \cite{MR715658} is written in the
following form: 
for any $f$ and its associated Maxwellian state $M$ with same mass,
momentum and temperature
$$ 
{\cal D}(f) \ge \lambda \, \rho \, \bigg[ {\cal H}(f) - {\cal H}(M) \bigg],
$$
where $\rho$ is the density (mass of $f$) and $\lambda >0$ is a
``{\it suitable constant}''. 

We shall now develop this very general statement into a layer of more
specified conjectures. Let us fix $\rho=1$ without loss of generality.

In the case when the constant $\lambda$ only depends on the collision
kernel $B$, the temperature of $M$ (or $f$), and some bound on the entropy of
$f$ ({\it i.e.}, only the basic physical {\it a priori} estimates), we
shall call this inequality the {\bf strong form of Cercignani's
  conjecture}.


In the case when the constant $\lambda$ also depends on some additional
bounds on $f$ (typically of smoothness, moments and lower bounds), we shall
call such an inequality the {\bf weak form of Cercignani's conjecture}. Let
us point out that it is of crucial importance to know whether the bounds
used can be shown to be propagated by the Boltzmann equation, in order to
be able to ``apply'' the weak form of Cercignani's conjecture to the
relaxation to equilibrium of its solutions. This of course guides which
bounds are natural or not.

In the slightly different case when the following inequality holds
$$ 
{\cal D}(f) \ge \lambda_\eps \, \big[ {\cal H}(f) - {\cal H}(M)
\big]^{1+\eps}, \quad \eps >0
$$
we shall speak of the {\bf $\eps$-polynomial Cercignani's conjecture} and
it can be divided again into weak and strong versions according to how much
the constant $\lambda_\eps$ depends on $f$.

Finally a strictly similar hierarchy of conjectures can be formulated on
the Landau entropy production functional, and we shall call it {\bf
  Cercignani's conjecture for the Landau equation}.

\subsection{A linearized counterpart to the conjecture}
\label{sec:line-count-conj}

A natural linearized counterpart of Cercignani's conjecture for the
Boltzmann or Landau equation consists in replacing the entropy production
functional and the Boltzmann entropy by their linearized approximation,
{\em i.e.}, respectively the Dirichlet form of the collision operators
discussed above and the $L^2(M)$ norm. This spectral gap question was
already well-known for a long time and used by Cercignani as an inspiration
and supportive argument for his conjecture in \cite{MR715658}.
So let us call this the {\bf linearized Cercignani's conjecture}:
$$ 
D(h) \ge \lambda \, \| h - \Pi(h)\|_{L^2(M)} ^2,
$$
where $\Pi$ denotes the orthogonal projector in $L^2(M)$ onto the null
space of the linearized collision operator, and $\lambda$ only
depends on the collision kernel $B$ and the temperature of $M$. 

Note that due the linear homogeneity of this relation, no weak version 
(with constant depending on the function $h$) would make sense.

Again obviously the same question can be asked on the Dirichlet form of the
Landau collision operators and leads to the {\bf linearized Cercignani's
  conjecture for the Landau equation}.

\subsection{Comparison with differential operators}
\label{sec:comp-with-diff}

In the light of the comparison we have made with usual differential
operators, a functional inequality interpretation of Cercignani's
conjecture is the following. Its nonlinear form is an intricate (because of
strong nonlinearity and average over additional angular variables)
amplified version of a logarithmic Sobolev inequality. Its linearized form
is an intricate (because again of average over additional angular
variables) amplified version of a Poincar\'e inequality.


\subsection{Use of the collision kernel, 
weighted forms of the conjecture}
\label{sec:use-collision-kernel}

A natural generalization of the conjecture is to know whether an
inequality of the form
$$ 
{\cal D}(f) \ge \lambda \, \bigg[ {\cal H}_w(f) - {\cal H}_w(M)
\bigg]
$$
holds, where 
\[
{\cal H}_w(f) := \int f \log f \, w(v) \, dv
\]
for some weight function $w>0$. This requires the additional moment
normalization 
\[
\int f \, w \, |v|^2 \, dv = \int M \, w \, |v|^2 \, dv
\]
in order for the relative entropy to satisfy
\[
{\cal H}_w (f) - {\cal H}_w(M) \ge 0.
\]
We shall call such an inequality a {\bf weighted Cercignani's
  conjecture}.

In the linearized case, a natural conjecture is similarly
$$ 
D(h) \ge \lambda \, \| (h - \Pi(h)) \, w \|_{L^2(M)} ^2,
$$
for some weight function $w$, and we shall call such an inequality a {\bf
  linearized weighted Cercignani's conjecture}.

\subsection{Semigroup form of the conjecture}
\label{sec:semigr-form-conj}

Another natural related question is the following. Cercignani's conjecture
at the end of the day is a purely functional inequality, and has nothing to
do with the solutions of the Boltzmann equation. However the main
application of this conjecture is of course the rate of convergence to
equilibrium for the Boltzmann equation. Hence a natural question is whether
$$ 
\forall \, t \ge 0, \quad {\cal D}(f_t) \ge \lambda \, \bigg[ {\cal H}(f_t) -
{\cal H}(M) \bigg]
$$
for any solution $(f_t)_{t \ge 0}$, or for a subset of solutions to
the spatially homogeneous Boltzmann equation. We shall call this the
{\bf strong semigroup Cercignani conjecture}.

This conjecture is of course related to the weak functional form of the
conjecture in the sense that if appropriate conditions are found on $f$ for
which the latter holds, and these conditions are shown to be propagated by
the solutions to the Boltzmann equation, then this strong semigroup form of
the conjecture will hold. However it is possible to imagine weird
conditions created by the semigroup that cannot be simply identified at the
purely functional level.

Let us first remark that Cercignani's conjecture does not exhaust the issue
of the rate of relaxation towards equilibrium, since it is possible that
an exponential convergence occurs for the semigroup even without any
functional inequality\footnote{Think for instance of the trivial example of
  a nonsymmetric matrix $A$ with only negative eigenvalues: its Dirichlet
  form will not control anything from below since it has no sign, although
  the semigroup $e^{tA}$ will obviously relax exponentially fast towards
  zero.}. Hence it is natural to weaken again the semigroup conjecture in
the following form. The {\bf weak form of the semigroup Cercignani
  conjecture} is the following: is it true that
$$ 
\forall \, t \ge 0, \quad \big| {\cal H}(f_t) -
{\cal H}(M) \big| \le C \, e^{-C' \, t}
$$
for some positive constants $C, C'>0$, and any solutions $(f_t)_{t \ge
  0}$, or for a subset of solutions to the spatially homogeneous
Boltzmann equation.

Let us add as a final remark that obvious extension to these semigroup
forms of the conjecture can be drawn for the spatially homogeneous Landau
equation.

\subsection{Relation with mean-field limit}
\label{sec:relation-with-mean}

Let us mention the different but closely related question raised by Kac
\cite{MR0084985}. In this paper Kac introduced a many-particle jump process
whose mean-field limit is the spatially homogeneous Boltzmann equation (he
also introduced the mathematical formalization of the by-now well-known
notion of propagation of chaos for a many-particles system). One of the
goals of such a derivation was the understanding of the asymptotic behavior
of the nonlinear Boltzmann equation through the linear many-particles
system. Hence it is natural to search for a Cercignani's conjecture {\bf at
  the level of the many-particles jump process}.  Even if this process is
linear the conjecture can be searched in nonlinear form or linear form, the
difference being the kind of functional which is used for measuring the
spectral gap (the Boltzmann entropy or an $L^2$ norm). The specific new
difficulty is to track the dependency of the estimates in terms of the
number of particles (since the final goal is originally to pass to the
limit).  For recent results in the ``linear'' case (spectral gap in $L^2$)
see for instance \cite{MR1743614,MR1825150,MR2020418,MR1961868} which
solves the problem (note however that since the $L^2$ norm does not
``tensorize'' correctly in high dimension, it is not possible to pass to
the limit in these estimates), and see
\cite{MR1964379,MR2580955,MISCHLER:2010:HAL-00447988:1} for some 
progress in the ``nonlinear'' case.

\section{Negative results at the functional level}

\subsection{Counterexample in the Maxwell molecules case 
with only mass, energy and entropy bounds}
\label{sec:counterexample-with}

The explicit solutions constructed by Bobylev
\cite{MR762269,MR1128328} show in particular that the exponential rate
of relaxation can be arbitrarily slow if one assumes only finite mass,
energy and entropy (in the cutoff Maxwell molecules case $B=1$). The
question was then: would the conjecture be true with better growth of
the collision kernel than $B=1$, for instance in the important
physical case of hard spheres $B=|v-v_*|$?

\subsection{Counterexample in the hard spheres case with higher moment
  bounds}
\label{sec:counterexample-hard}

Then Wennberg \cite{MR1450762} proved in 1997 that the conjecture was
false also in the case of hard spheres. However the counterexample of
Wennberg does not have infinitely many moments. And since it was known
since \cite{MR1233644} that hard potentials with cutoff ($\gamma>0$
and no angular singularity) would produce such higher moments at every
order for any positive time, the next natural question was whether
such counterexamples would hold with infinitely many moments.

\subsection{Counterexample with infinitely many moment bounds}
\label{sec:counterexample-with-1}

It was finally shown by Bobylev and Cercignani in \cite{MR1675366}
that, for hard potentials with angular cutoff and with $\gamma <2$,
the estimate
\[ 
{\cal D}(f) \ge C(f) \, \Bigl[ {\cal H}(f) - {\cal H}(M) \Bigr] 
\] 
does not hold uniformly over any class of distributions $f$ which
have fixed mass, momentum and energy, are bounded below by a uniform Maxwellian distribution,
and satisfy
\[
\| f \|_{H^k} \le S_k, \quad  \|f\|_{L^1(1+|v|^s)} \le M_s
\]
for some given sequences $M_s>0$, $s \in \N$ and $S_k>0$, $k \in \N$.
This is still true if an arbitrary number of the
constants $M_s$ and $S_k$ are chosen arbitrarily close to their
equilibrium value ({\it i.e.}, the value of the corresponding moment or
Sobolev norm for the Maxwellian distribution associated with $f$).

\section{Positive results at the functional level}

\subsection{Early attempts}

A first result was proved in 1979 by Aizenmann and Bak in~\cite{AB} on a
nonlinear model bearing some reminiscence to the Boltzmann operator,
namely a coagulation-fragmentation kernel with constant rate. The proof
uses two very cleverly devised convexity inequalities and this is probably
the first example of a linear relation between entropy and entropy
production for a non-diffusive Boltzmann-like equation.
 
For the Boltzmann and Landau equations, an inequality linking the entropy
dissipation towards an $L^1$-type distance to the space of Maxwellian
(instead of the entropy!) was established by the first author
in~\cite{desvilleCMP}, for densities bounded below by a Maxwellian.

A decisive improvement was the discovery of a link between the entropy and
the entropy dissipation of the Boltzmann equation, due to Carlen and
Carvalho in the beginning of the 90s (Cf. \cite{CC1, CC2}): this link was
however still far from linear, the relation was of the form
$$
{\cal D}(f) \ge \Theta\bigl( {\cal H}(f)  - {\cal H}(M) \bigr), 
$$
for some nonlinear function $\Theta$, but no information was given on the
behavior of the function $\Theta$ at zero, and it turns out that its
behavior was much worse than linear. The inequality was proved in the case
of pseudo-Maxwell molecules and its proof made crucial use of the theory of
Wild sums. As an elaborated consequence of this breakthrough let us also
mention the paper \cite{MR1475461} where the authors give (for a constant
collision kernel) the first proof of the convergence of the entropy towards
its equilibrium value with minimal assumptions on the initial datum
(finite mass, energy and entropy).

\subsection{Landau equation with (over)-Maxwellian molecules}

The first positive result obtained with a linear dependence between
entropy and entropy dissipation (for Landau or Boltzmann operators) is
due to the first and third author of this paper, and concerns Landau's
operator
in the case {when} {$\Phi(|v-v_*|) \ge C_\Phi >0$}
(sometimes called ``over-Maxwellian''): 
$$ 
{\mathcal{D}}(f) \ge \lambda\, \Big[ {\mathcal{H}}(f) -
{\mathcal{H}}(M) \Big],
$$ 
with {$\lambda>0$} depending on $C_\Phi$ and the dimension $d$ only (for
all {$f$} with fixed mass, momentum, energy and upper bound on the
entropy). In other words, the strong conjecture holds for Landau equation
with over-Maxwellian molecules.

This result is obtained by observing that (for some $\lambda'>0$) { $$
  {\mathcal{D}}(f) \ge \lambda'\, {\mathcal{D}}_{FP}(f) , $$} and then by
using the logarithmic Sobolev inequality.  

The relationship between ${\mathcal{D}}(f)$ and ${\mathcal{D}}_{FP}(f)$ can
be obtained (Cf. \cite{MR1737548}) either thanks to an almost explicit
computation in which spherical coordinates are used, or thanks to a variant
of the estimates of \cite{desvilleCMP} (in which differential operators
have been replaced by integral operators).

 \subsection{Boltzmann equation: ``almost linear'' relation}

 Successive results by Toscani and the third author of this survey
 (Cf. \cite{TV1, TV2, MR1964379}) then led to a result valid for the
 Boltzmann equation with any reasonable collision kernel. That is,
 when {$B \ge C \min \{ |v-v_*|^\gamma, |v-v_*|^{-\beta} \} >0$, }
 with {$\beta, \gamma\ge 0$} (with or without cutoff, that is for any
 $\nu \in (-\infty,2)$) the following ``almost linear link'' exists
 between entropy and entropy production:
\[  
{\mathcal{D}}(f) \ge  \lambda_{\var}(f)
 \, \Big[ {\mathcal{H}}(f) - {\mathcal{H}}(M) \Big]^{1+\var}, 
\]
where $\lambda_{\var}(f)$ (possibly going to $0$ as $\eps \to 0$)
depends on {$\beta, \gamma, C,$} $\| f \|_{H^{s(\var)}}$, $ \| f
(1+|v|)^{k(\var)} \|_{L^1}$ for some $s(\eps)$ and $k(\eps)$ depending on
$\eps$ (and blowing up as $\eps \to 0$), and some Gaussian (or better)
lower bounds.  In other terms, the (weak) $\var$-polynomial
Cercignani's conjecture holds for the Boltzmann operator with all
reasonable collision kernels.

 \subsection{Boltzmann equation: linear relation}

 Finally, in the special (and unfortunately non physical) case when
 $\Phi(|v-v_*|) \ge C_\Phi \, (1 + |v-v_*|^2) >0$, the original
 (strong) Cercignani's conjecture was proved to hold by the third
 author of this paper (Cf. \cite{MR1964379}), that is:
$$  
{\mathcal{D}}(f) \ge \lambda \, \Big[ {\mathcal{H}}(f) -
{\mathcal{H}}(M) \Big],
$$
with {$\lambda>0$} only depending on $C_\Phi$ and the energy and entropy of
$f$. 


\subsection{Spectral gaps of linearized operators}

The origin of the study of the spectral properties of the linearized
Boltzmann collision operator can be traced back all the way to Hilbert
\cite{Hilb:EB:12} in his work on the integral operators. The integral
collision operator of Boltzmann was a key example and motivation from
physics there, and Hilbert introduced (in the case of hard spheres, {\it
  i.e.}, $\gamma=1$ and $\nu <0$ in dimension $d=3$) the splitting between
the local and non-local parts of the linearized operator, and proved the
``complete continuity'' (compactness in today's words) of its non-local
part. This result was important in the construction, together with the
Fredholm theory, of what is now called the ``Hilbert expansions''.

The second key step is due to Carleman \cite{Carl:fond:57}. In this book he
introduced (still in the hard spheres case) the use of Weyl's theorem
(stating that under certain assumptions, the essential spectrum is
unchanged by a compact perturbation of the operator), in order to prove the
existence of a spectral gap of the linearized collision operator. Grad
\cite{Grad63} then generalized this result to a broader class of collision
kernels (the so-called hard potentials with cutoff, $\gamma\in (0,1]$, $\nu
<0$). Further generalizations of this non constructive approach were
provided by Caflisch \cite{Cafl:LBE1:80} and Golse and Poupaud
\cite{GolPo:Cras:86} who proved that in the case of soft potentials
with angular cutoff ($\gamma \in (-d,0)$ and $\nu <0$) there is no spectral
gap but still a weighted form of the conjecture holds (in weaker norms than
the ambiant norm).

Another direction of research was investigated by  Wang-Chang and Uhlenbeck
\cite{WCUh:LBE:70} and Bobylev \cite{MR1128328}, who obtained a complete
diagonalization of the linearized collision operator in the Maxwell
molecules case ($\gamma=0$). 

In the paper \cite{Baranger-Mouhot-2005}, by Baranger and the second author
of this survey, the first constructive estimates
were given for the spectral gap of the hard spheres model, by introducing a
geometric method and the idea of using ``intermediate'' collisions in the
regions where the collision kernel vanishes, and then rely on the explicit
diagonalization of the Maxwell case above (for a constant collision
kernel). The paper \cite{MR2254617} then developed a whole series of
constructive coercivity estimates, most of them sharp, and gives a full
answer to the linearized Cercignani's conjecture in the cutoff case (and
some ellipticity estimate in the non-cutoff case, see below).

The interested reader can find more details in the review paper
\cite{MR2301289}.

\subsection{Elliptic estimates in the long-range case}
\label{sec:ellipt-estim-long}

A different but related subject is the issue of proving ellipticity
estimate on the entropy production functional in the long-range case
$\nu \in (0,2)$. The first results are due to Lions \cite{MR1649477}
and then to the first and third authors of this paper, together with
Alexandre and Wennberg (\cite{MR1715411} and \cite{ADVW:00}) and can
be roughly stated as the following inequality
\[ 
{\cal D}(f) \ge C \, \Bigl( \| f \|_{H^{\nu/2} _{loc}} - \left\| f (1+|v|^s)
\right\|_{L^1} \Bigr)
\]
for some $s \in \N$ and some constant $C>0$, and where the exponent $\nu/2$
of regularization is related to the singularity of the collision kernel.
These regularity estimates were used in a crucial way in the papers
\cite{AlVi:02,MR2037247} by Alexandre and the third author of this
survey in the Cauchy theory of the long-ranged Boltzmann equation and
its grazing collision limit towards the Landau equation. Earlier an
estimate based on the entropy production had already been used in a
crucial way in the paper \cite{MR1650006} in the Cauchy theory of the
spatially homogeneous long-ranged Boltzmann with soft potentials.

The corresponding (and sharper) regularity estimates for the Landau
equation were derived in the papers \cite{MR1737547,MR1737548} by the
first and third authors of this survey.

Now if we consider the linearized case, the first related result for
the Boltzmann collision operator is due to Pao
\cite{Pao:TSnoncutoff:74} who proved the compactness of the resolvent
in the case of inverse power-law interaction force $F(z) = z^{-s}$
with $s > 3$ in dimension $d=3$ (this paper was hard to read because
of the use of pseudo-differential calculus, and somewhat controversial,
see \cite{Klaus:TSnoncutoff:77}). However, some more general and
constructive versions of these results were recovered by the new
estimates proven by Strain and the second author of this paper in
\cite{MR2322149}, confirming that the latter paper was fully
correct. In short the paper \cite{MR2322149} proves that a ``gain'' of
$\nu$ polynomial weight occurs in the long-ranged case due to the
ellipticity of the operator: roughly speaking the local regularization
improves the behavior for large velocities. The recent work by
Gressman and Strain \cite{MR2629879} finally gives a sharp (although
intricate) characterization of the elliptic norm associated with the
Dirichlet form, and answers a conjecture raised in \cite{MR2322149}
about the spectral gap in the non-cutoff case (see below).

For the linearized Landau operator, the non-constructive approach
based on Weyl's theorem was used successfully by Degond and Lemou
\cite{DeLe:LLE:97} to prove an elliptic coercivity estimate. The sharp
(non-isotropic) norm was then devised by Guo in \cite{Guo-2002-II} and
some constructive versions of this coercivity result were given in
\cite{MR2322149} (see also \cite{herau-pravda} for related sharp
hypoellipticity results).

\section{Results on the semigroup Cercignani conjecture}
\label{sec:count-exampl-semigr}


Since up to this date no results are known on the strong semigroup
Cercignani's conjecture, we shall review the results on the weak
semigroup conjecture (about the rate of decay of the semigroup).

The first paper which tackled the issue of exponential relaxation to
equilibrium for data not necessarily close to equilibrium is
\cite{MR900501} where a non constructive result is given by Arkeryd,
Esposito and Pulvirenti in the hard spheres case for general solutions
in weighted $L^\infty$ spaces in the spatially homogeneous case, or
solutions weakly inhomogeneous (close to the latter setting). This
approach was generalized to an $L^p$ setting in \cite{Wenn:stab:93} by
Wennberg.

In the papers \cite{MR1991033,MR2546739}, Carlen, Carvalho and Lu give
sharp lower (and upper) bounds on the rate of relaxation to
equilibrium for the semigroup of the spatially homogeneous Boltzmann
equation in the case of Maxwell molecules and soft potentials. Let us
also mention the related important works by Carlen, Carvalho and
Gabetta \cite{MR1725612,MR2119283}, more focused on the specific issue
of the Wild sums.

The paper \cite{MR1669689} by Carlen, Gabetta and Toscani then gave a
proof of exponential convergence towards equilibrium in the particular
case of Maxwell molecules with cutoff, assuming moments and regularity
bounds on the initial datum.

In the Maxwell molecules case with angular cutoff ($\gamma=0$, $\nu
<0$ in the classification above), Carlen and Lu proved in
\cite{MR1991033} that when only $2$ moments are assumed on the initial
datum, then a rate of convergence algebraic and as slow as wanted can
be explicitely obtained, and it depends very precisely on the
integrability of the initial datum, in the sense that it depends on
the function $\varphi=\varphi(v) >0$ (going to infinity as $v$ goes to
infinity) such that
\[
\int f_{in} \, (1+|v|^2) \, \varphi \, dv < +\infty. 
\]
This was the first example of a solution to the Boltzmann equation that
{\em actually relaxes slower than exponentially towards equilibrium}.

In the soft potential case with or without angular cutoff ($\gamma \in
(-d,0)$ and $\nu <2$), Carlen, Carvalho and Lu \cite{MR2546739} proved (among
other things) that for solutions with only $2+|\gamma|$ moments then again
the convergence (in $L^1$) can be algebraic and arbitrarily slow.

More recently a whole program of research has been carried out by the
first and third author in order to develop a theory of relaxation to
equilibrium in the large for inhomogeneous kinetic equations, and in
particular for the full Boltzmann equation in confined domain under
{\it a priori} smoothness and moments assumptions (see in particular
the papers
\cite{Desvillettes-Villani-2001,Desvillettes-Villani-2005}). The rate
is not exponential, but still almost exponential (in the sense of a
polynomial convergence for any polynomial). These works have provided
new insights on how transport effects combined with the thermalization
process in velocity can yield convergence towards the global
equilibrium and they have given birth to the new {\em hypocoercivity}
theory.

Finally the last episode in this long research line on the rate of
decay in the Boltzmann $H$ theorem is the development of new tools in
spectral theory for non-symmetric operators in order to systematically
enlarge the functional space of the decay on a linearized semigroup
(this is used for connecting the linearized theory
\cite{Baranger-Mouhot-2005,Mouhot-Neumann} in small $L^2(M^{-1})$-type
spaces with the nonlinear theory
\cite{MR1964379,Desvillettes-Villani-2005} in $L^1$ type spaces). A
``final answer'' (constructive exponential relaxation of the relative
entropy under reasonable {\it a priori} estimates) was given in
\cite{MR2197542} for hard spheres in the spatially homogeneous case,
and then in the preprint by Gualdani, Mischler and the second author
of this survey \cite{GUALDANI:2010:HAL-00495786:1} for the
inhomogeneous case (confined in the torus).

\section{New conjectures about the conjecture}
\label{sec:conjectures}

We shall now present a few conjectures 
related to the results presented in this paper, some of them already formulated,
some of them new. The work in progress \cite{DM-inprogress} aims at giving 
partial answers to some of these conjectures.

\subsection{The use of grazing collisions}

In \cite{MR1964379} the following conjecture is formulated:
\smallskip

\noindent {\bf Conjecture~1 (Villani).} {\it For a collision kernel with
  growth of order $\gamma \in (-d,+\infty)$ and singularity of order $\nu
  \in [0,2]$ (where $\nu=2$ formally plays the role of the Landau collision
  operator), the strong form of Cercignani's conjecture is true if and only
  if $\gamma + \nu \ge 2$.}
\smallskip

This conjecture formally extends  the result proved in
\cite{MR1964379} for ``superquadratic" collision kernels (formally $\nu=0$
and $\gamma \ge 2$), as well as the result obtained in \cite{MR1737548} 
for the Landau collision operator with hard
potential (formally $\nu=2$ and $\gamma \ge 0$).

In \cite{MR2322149} it is conjectured the following 
\smallskip

\noindent {\bf Conjecture~2 (Mouhot-Strain).} {\it For a collision kernel with
  growth of order $\gamma \in (-d,+\infty)$ and singularity of order $\nu
  \in [0,2]$ (where $\nu=2$ formally plays the role of the Landau collision
  operator), the strong form of linearized Cercignani's conjecture
  (existence of a spectral gap for the linearized operator) is true if and
  only if $\gamma + \nu \ge 0$.}
\smallskip

The direct implication in this conjecture was proved in
\cite{MR2322149}, and recently the work \cite{MR2629879} (and the
preprints related to this annoucement note) has solved this conjecture by
providing a sharp characterization of the norm associated with the
Dirichlet form of the operator. Let us also mention as an example of the
fertility of Cercignani's conjecture that this conjecture also has inspired
the work \cite{MRS} about fractional Poincar\'e inequalities.

An interesting open question which now calls for further understanding is
why there is such a ``gap'' between the condition $\gamma+\nu\ge 2$ for a
``nonlinear'' spectral gap, and the condition $\gamma + \nu \ge 0$ for a
linearized spectral gap (this is somehow reminiscent of the gap between the
conditions for a logarithmic Sobolev inequality to hold or a Poincar\'e
inequality to hold).

\subsection{Use of (general) collision kernels}

In view of the preceding conjectures about the influence of the parameters
$\gamma$ and $\nu$ (both at the nonlinear and linearized levels), it now
seems natural to conjecture the following generalized weighted form of
Cercignani's conjecture: 
\smallskip

\noindent
{\bf Conjecture 3.}  {\it For a collision kernel with
  growth of order $\gamma \in (-d,+\infty)$ and singularity of order $\nu
  \in [0,2]$ (where $\nu=2$ formally plays the role of the Landau collision
  operator), one has 
  \[ {\cal D}(f) \ge \lambda \, \int f \log \frac{f}{M} \,
  (1+|v|)^{\gamma+\nu-2} \, dv
\]
at the nonlinear level and 
\[
D(h) \ge \lambda' \, \| (h -\Pi h) (1+|v|)^{\gamma+\nu} \|_{L^2(M)}
\]
at the linearized level (where $h = (f-M)/M$ denotes the rescaled
fluctuation around the equilibrium $M$, and $\Pi$ the $L^2(M)$-orthogonal
projection onto the collision invariants $1,v_1, \dots ,v_d,|v|^2$).}

\subsection{Use of the tail of the distribution}
 
In view of the condition for a measure $\mu=e^{-V}$ to satisfy a Poincar\'e
inequality (essentially that $V$ grows faster than linearly at infinity) or
a logarithmic Sobolev inequality (essentially that $V$ grows faster than
quadratically at infinity), see for instance \cite{Vill:tspt}, it is
natural to ask whether some bounds on the tail of the distribution with
faster decay than polynomial could help in order to obtain Cercignani's
conjecture (this would be compatible with the counterexamples of
\cite{MR1675366}):
\smallskip

\noindent
{\bf Conjecture 4.}  {\it For a collision kernel with
  growth of order $\gamma \in (-d,+\infty)$ and singularity of order $\nu
  \in [0,2]$ (where $\nu=2$ formally plays the role of the Landau collision
  operator), and for a distribution $f$ with exponential tail 
\[
f(v) \sim C \, e^{-a |v|^\eta}, \quad v \sim \infty
\]
for some constants $a, C>0$ and $\eta \in (0,2]$,  one has 
  \[ {\cal D}(f) \ge \lambda \, \int f \log \frac{f}{M} \,
  (1+|v|)^{\gamma+\nu+\eta-2} \, dv
\]
at the nonlinear level.}
\smallskip

Hence if this conjecture were true (and taking for granted Villani's
conjecture), it would mean that Cercignani's conjecture would hold for
hard spheres and inverse power-laws interactions in dimension $d=3$
for initial data with strong enough exponential tail (a Gaussian tail
would always work). These exponential tails are known to
be propagated by the spatially homogeneous nonlinear Boltzmann
equation for hard spheres and hard potentials with cutoff
\cite{MR1478067,MR2197542} (see also \cite{MR2533928} in the case of
Gaussian tails).

Concerning general initial data, one has to turn to the theory of
appearance of exponential moments, as developed in
\cite{MR2197542,MR2264623}. However it only provides at the moment an
appearance of an $L^1$ decay of the form $e^{-a |v|^\eta}$ with
$\eta=\gamma/2$. Even with the improvement to $\eta=\gamma$ in the recent
work in progress \cite{canizo-mouhot}, note that the hard spheres case (say
$\gamma=0$ and $\nu=0$) is exactly borderline for Cercignani's conjecture,
and for inverse power-laws in dimension $d=3$ where $\gamma = (s-5)/(s-1)$
and $\nu=2/(s-1)$, then $2 \gamma + \nu = 2 - 6 / (s-1)$ is always less
than $2$, and hence this would not be sufficient for Cercignani's
conjecture.

Of course a natural further question would be to know whether the grazing
collisions $\nu >0$ can help for the appearance of exponential moments (for
instance with $\eta=\gamma + \nu$ instead of $\eta=\gamma$)\dots As the
reader sees, the story opened by Cercignani's conjecture is far from
finished, and it is likely that many more nice surprises are to come.

\bigskip
\noindent {\bf{Acknowledgments:}} The authors wish to thank the ANR
grant CBDif for support. The second author wishes to thank the Award
No. KUK-I1-007-43, funded by the King Abdullah University of Science
and Technology (KAUST) for the funding provided for his repeated
visits at the University of Cambridge during the autumn 2009 and the
spring 2010.\smallskip




\bibliographystyle{acm} 
\bibliography{DM-EEP}

\def\cprime{$'$} \def\cprime{$'$}
\begin{thebibliography}{10}

\bibitem{AB}
{\sc Aizenman, M., and Bak, T.}
\newblock Convergence to equilibrium in a system of reacting polymers.
\newblock {\em Comm. Math. Phys. 65\/} (1979), 203--230.

\bibitem{ADVW:00}
{\sc Alexandre, R., Desvillettes, L., Villani, C., and Wennberg, B.}
\newblock Entropy dissipation and long-range interactions.
\newblock {\em Arch. Rational Mech. Anal. 152}, 4 (2000), 327--355.

\bibitem{AlVi:02}
{\sc Alexandre, R., and Villani, C.}
\newblock On the {B}oltzmann equation for long-range interactions.
\newblock {\em Comm. Pure Appl. Math. 55}, 1 (2002), 30--70.

\bibitem{MR2037247}
{\sc Alexandre, R., and Villani, C.}
\newblock On the {L}andau approximation in plasma physics.
\newblock {\em Ann. Inst. H. Poincar\'e Anal. Non Lin\'eaire 21}, 1 (2004),
  61--95.

\bibitem{canizo-mouhot}
{\sc Alonso, R., Ca\~nizo, J., Gamba, I., and Mouhot, C.}
\newblock Exponential moments for the spatially homogeneous {B}oltzmann
  equation.
\newblock Work in progress.

\bibitem{MR900501}
{\sc Arkeryd, L., Esposito, R., and Pulvirenti, M.}
\newblock The {B}oltzmann equation for weakly inhomogeneous data.
\newblock {\em Comm. Math. Phys. 111}, 3 (1987), 393--407.

\bibitem{MR1055522}
{\sc Arsen{\cprime}ev, A.~A., and Buryak, O.~E.}
\newblock On a connection between the solution of the {B}oltzmann equation and
  the solution of the {L}andau-{F}okker-{P}lanck equation.
\newblock {\em Mat. Sb. 181}, 4 (1990), 435--446.

\bibitem{BK}
{\sc Bakry, D., and Emery, M.}
\newblock Diffusions hypercontractives.
\newblock {\em S\'em. Proba. XIX, Lecture Notes in Math. 1123\/} (1985),
  177--206.

\bibitem{Baranger-Mouhot-2005}
{\sc Baranger, C., and Mouhot, C.}
\newblock Explicit spectral gap estimates for the linearized {B}oltzmann and
  {L}andau operators with hard potentials.
\newblock {\em Rev. Mat. Iberoamericana 21}, 3 (2005), 819--841.

\bibitem{MR762269}
{\sc Bobylev, A.~V.}
\newblock Exact solutions of the nonlinear {B}oltzmann equation and the theory
  of relaxation of a {M}axwell gas.
\newblock {\em Teoret. Mat. Fiz. 60}, 2 (1984), 280--310.

\bibitem{MR1128328}
{\sc Bobylev, A.~V.}
\newblock The theory of the nonlinear spatially uniform {B}oltzmann equation
  for {M}axwell molecules.
\newblock In {\em Mathematical physics reviews, {V}ol.\ 7}, vol.~7 of {\em
  Soviet Sci. Rev. Sect. C Math. Phys. Rev.} Harwood Academic Publ., Chur,
  1988, pp.~111--233.

\bibitem{MR1478067}
{\sc Bobylev, A.~V.}
\newblock Moment inequalities for the {B}oltzmann equation and applications to
  spatially homogeneous problems.
\newblock {\em J. Stat. Phys. 88}, 5-6 (1997), 1183--1214.

\bibitem{MR1675366}
{\sc Bobylev, A.~V., and Cercignani, C.}
\newblock On the rate of entropy production for the {B}oltzmann equation.
\newblock {\em J. Stat. Phys. 94}, 3-4 (1999), 603--618.

\bibitem{boltzmann}
{\sc Boltzmann, L.}
\newblock Weitere studien uber das w\"arme gleichgenicht unfer gasmol\"akuler.
\newblock {\em Sitzungsberichte der Akademie der Wissenschaften 66\/} (1872),
  265--370.
\newblock Translation: Further studies on the thermal equilibrium of gas
  molecules, in Kinetic Theory 2, 88–174, Ed. S.G. Brush, Pergamon, Oxford
  (1966).

\bibitem{Cafl:LBE1:80}
{\sc Caflisch, R.~E.}
\newblock The {B}oltzmann equation with a soft potential. {I}. {L}inear,
  spatially-homogeneous.
\newblock {\em Comm. Math. Phys. 74}, 1 (1980), 71--95.

\bibitem{Carl:fond:57}
{\sc Carleman, T.}
\newblock {\em Probl\`emes Math\'ematiques dans la Th\'eorie Cin\'etique des
  Gaz}.
\newblock Almqvist \& Wiksell, 1957.

\bibitem{CC1}
{\sc Carlen, E.~A., and Carvalho, M.}
\newblock Strict entropy production bounds and stability of the rate of
  convergence to equilibrium for the {B}oltzmann equation.
\newblock {\em J. Stat. Phys. 67}, 3-4 (1992), 575--608.

\bibitem{CC2}
{\sc Carlen, E.~A., and Carvalho, M.}
\newblock Entropy production estimates for {B}oltzmann equations with
  physically realistic collision kernels.
\newblock {\em J. Stat. Phys. 74}, 3-4 (1994), 743--782.

\bibitem{MR1725612}
{\sc Carlen, E.~A., Carvalho, M.~C., and Gabetta, E.}
\newblock Central limit theorem for {M}axwellian molecules and truncation of
  the {W}ild expansion.
\newblock {\em Comm. Pure Appl. Math. 53}, 3 (2000), 370--397.

\bibitem{MR2119283}
{\sc Carlen, E.~A., Carvalho, M.~C., and Gabetta, E.}
\newblock On the relation between rates of relaxation and convergence of {W}ild
  sums for solutions of the {K}ac equation.
\newblock {\em J. Funct. Anal. 220}, 2 (2005), 362--387.

\bibitem{MR2580955}
{\sc Carlen, E.~A., Carvalho, M.~C., Le~Roux, J., Loss, M., and Villani, C.}
\newblock Entropy and chaos in the {K}ac model.
\newblock {\em Kinet. Relat. Models 3}, 1 (2010), 85--122.

\bibitem{MR2020418}
{\sc Carlen, E.~A., Carvalho, M.~C., and Loss, M.}
\newblock Determination of the spectral gap for {K}ac's master equation and
  related stochastic evolution.
\newblock {\em Acta Math. 191}, 1 (2003), 1--54.

\bibitem{MR2546739}
{\sc Carlen, E.~A., Carvalho, M.~C., and Lu, X.}
\newblock On strong convergence to equilibrium for the {B}oltzmann equation
  with soft potentials.
\newblock {\em J. Stat. Phys. 135}, 4 (2009), 681--736.

\bibitem{MR1475461}
{\sc Carlen, E.~A., Carvalho, M.~C., and Wennberg, B.}
\newblock Entropic convergence for solutions of the {B}oltzmann equation with
  general physical initial data.
\newblock {\em Transport Theory Statist. Phys. 26}, 3 (1997), 373--378.

\bibitem{MR1669689}
{\sc Carlen, E.~A., Gabetta, E., and Toscani, G.}
\newblock Propagation of smoothness and the rate of exponential convergence to
  equilibrium for a spatially homogeneous {M}axwellian gas.
\newblock {\em Comm. Math. Phys. 199}, 3 (1999), 521--546.

\bibitem{MR1991033}
{\sc Carlen, E.~A., and Lu, X.}
\newblock Fast and slow convergence to equilibrium for {M}axwellian molecules
  via {W}ild sums.
\newblock {\em J. Stat. Phys. 112}, 1-2 (2003), 59--134.

\bibitem{cercignani}
{\sc Cercignani, C.}
\newblock {\em Theory and application of the {B}oltzmann equation}.
\newblock Elsevier, New York, 1975.

\bibitem{MR715658}
{\sc Cercignani, C.}
\newblock {$H$}-theorem and trend to equilibrium in the kinetic theory of
  gases.
\newblock {\em Arch. Mech. (Arch. Mech. Stos.) 34}, 3 (1982), 231--241 (1983).

\bibitem{CZ1}
{\sc Csiszar, I.}
\newblock Information-type measures of difference of probability distributions
  and indirect observations.
\newblock {\em Stud. Sci. Math. Hung. 2\/} (1967), 299--318.

\bibitem{DeLe:LLE:97}
{\sc Degond, P., and Lemou, M.}
\newblock Dispersion relations for the linearized {F}okker-{P}lanck equation.
\newblock {\em Arch. Rational Mech. Anal. 138}, 2 (1997), 137--167.

\bibitem{MR1167768}
{\sc Degond, P., and Lucquin-Desreux, B.}
\newblock The {F}okker-{P}lanck asymptotics of the {B}oltzmann collision
  operator in the {C}oulomb case.
\newblock {\em Math. Models Methods Appl. Sci. 2}, 2 (1992), 167--182.

\bibitem{desvilleCMP}
{\sc Desvillettes, L.}
\newblock Entropy dissipation rate and convergence in kinetic equations.
\newblock {\em Comm. Math. Phys. 123}, 4 (1989), 687--702.

\bibitem{MR1165528}
{\sc Desvillettes, L.}
\newblock On asymptotics of the {B}oltzmann equation when the collisions become
  grazing.
\newblock {\em Transport Theory Statist. Phys. 21}, 3 (1992), 259--276.

\bibitem{MR1233644}
{\sc Desvillettes, L.}
\newblock Some applications of the method of moments for the homogeneous
  {B}oltzmann and {K}ac equations.
\newblock {\em Arch. Rational Mech. Anal. 123}, 4 (1993), 387--404.

\bibitem{DM-inprogress}
{\sc Desvillettes, L., and Mouhot, C.}
\newblock About {C}ercignani's conjecture.
\newblock Work in progress.

\bibitem{MR1737547}
{\sc Desvillettes, L., and Villani, C.}
\newblock On the spatially homogeneous {L}andau equation for hard potentials.
  {I}. {E}xistence, uniqueness and smoothness.
\newblock {\em Comm. Partial Differential Equations 25}, 1-2 (2000), 179--259.

\bibitem{MR1737548}
{\sc Desvillettes, L., and Villani, C.}
\newblock On the spatially homogeneous {L}andau equation for hard potentials.
  {II}. {$H$}-theorem and applications.
\newblock {\em Comm. Partial Differential Equations 25}, 1-2 (2000), 261--298.

\bibitem{Desvillettes-Villani-2001}
{\sc Desvillettes, L., and Villani, C.}
\newblock On the trend to global equilibrium in spatially inhomogeneous
  entropy-dissipating systems: the linear {F}okker-{P}lanck equation.
\newblock {\em Comm. Pure Appl. Math. 54}, 1 (2001), 1--42.

\bibitem{Desvillettes-Villani-2005}
{\sc Desvillettes, L., and Villani, C.}
\newblock On the trend to global equilibrium for spatially inhomogeneous
  kinetic systems: the {B}oltzmann equation.
\newblock {\em Invent. Math. 159}, 2 (2005), 245--316.

\bibitem{MR1743614}
{\sc Diaconis, P., and Saloff-Coste, L.}
\newblock Bounds for {K}ac's master equation.
\newblock {\em Comm. Math. Phys. 209}, 3 (2000), 729--755.

\bibitem{MR2533928}
{\sc Gamba, I.~M., Panferov, V., and Villani, C.}
\newblock Upper {M}axwellian bounds for the spatially homogeneous {B}oltzmann
  equation.
\newblock {\em Arch. Ration. Mech. Anal. 194}, 1 (2009), 253--282.

\bibitem{GolPo:Cras:86}
{\sc Golse, F., and Poupaud, F.}
\newblock Un r\'esultat de compacit\'e pour l'\'equation de {B}oltzmann avec
  potentiel mou. {A}pplication au probl\`eme de demi-espace.
\newblock {\em C. R. Acad. Sci. Paris S\'er. I Math. 303}, 12 (1986), 583--586.

\bibitem{Grad63}
{\sc Grad, H.}
\newblock Asymptotic theory of the {B}oltzmann equation. {II}.
\newblock In {\em Rarefied Gas Dynamics (Proc. 3rd Internat. Sympos., Palais de
  l'UNESCO, Paris, 1962), Vol. I}. Academic Press, New York, 1963, pp.~26--59.

\bibitem{MR2629879}
{\sc Gressman, P.~T., and Strain, R.~M.}
\newblock Global classical solutions of the {B}oltzmann equation with
  long-range interactions.
\newblock {\em Proc. Natl. Acad. Sci. USA 107}, 13 (2010), 5744--5749.

\bibitem{MR0372613}
{\sc Gross, L.}
\newblock Hypercontractivity and logarithmic {S}obolev inequalities for the
  {C}lifford {D}irichlet form.
\newblock {\em Duke Math. J. 42}, 3 (1975), 383--396.

\bibitem{GUALDANI:2010:HAL-00495786:1}
{\sc {G}ualdani, M.~P., {M}ischler, S., and {M}ouhot, C.}
\newblock {F}actorization for non-symmetric operators and exponential
  {H}-theorem.
\newblock ArXiv preprint 1006.5523 (2010).

\bibitem{Guo-2002-II}
{\sc Guo, Y.}
\newblock The {L}andau equation in a periodic box.
\newblock {\em Comm. Math. Phys. 231}, 3 (2002), 391--434.

\bibitem{herau-pravda}
{\sc H\'erau, F., and Pravda-Starov, K.}
\newblock Anisotropic hypoelliptic estimates for landau-type operators.
\newblock ArXiv preprint 1003.3265 (2010).

\bibitem{Hilb:EB:12}
{\sc Hilbert, D.}
\newblock Grundz{\"u}ge einer {A}llgemeinen {T}heorie der {L}inearen
  {I}ntegralgleichungen.
\newblock {\em Math. Ann. 72\/} (1912).
\newblock Chelsea {P}ubl., {N}ew {Y}ork, 1953.

\bibitem{MR1825150}
{\sc Janvresse, E.}
\newblock Spectral gap for {K}ac's model of {B}oltzmann equation.
\newblock {\em Ann. Probab. 29}, 1 (2001), 288--304.

\bibitem{MR0084985}
{\sc Kac, M.}
\newblock Foundations of kinetic theory.
\newblock In {\em Proceedings of the {T}hird {B}erkeley {S}ymposium on
  {M}athematical {S}tatistics and {P}robability, 1954--1955, vol. {III}\/}
  (Berkeley and Los Angeles, 1956), University of California Press,
  pp.~171--197.

\bibitem{Klaus:TSnoncutoff:77}
{\sc Klaus, M.}
\newblock Boltzmann collision operator without cut-off.
\newblock {\em Helv. Phys. Acta 50}, 6 (1977), 893--903.

\bibitem{CZ2}
{\sc Kullback, S.}
\newblock A lower bound for discrimination information in terms of variation.
\newblock {\em IEEE Trans. Inf. The. 4\/} (1967), 126--127.

\bibitem{landau-1936}
{\sc Landau, L.}
\newblock Die kinetische gleichung f\"ur den fall coulombscher wechselwirkung.
\newblock {\em Phys. Z. Sowjet. 154}, 10 (1936).
\newblock Translation: The transport equation in the case of Coulomb
  interactions, in D. ter Haar, ed., Collected papers of L.D. Landau, pp.
  163–170. Pergamon Press, Oxford, 1981.

\bibitem{MR1649477}
{\sc Lions, P.-L.}
\newblock R\'egularit\'e et compacit\'e pour des noyaux de collision de
  {B}oltzmann sans troncature angulaire.
\newblock {\em C. R. Acad. Sci. Paris S\'er. I Math. 326}, 1 (1998), 37--41.

\bibitem{MR1961868}
{\sc Maslen, D.~K.}
\newblock The eigenvalues of {K}ac's master equation.
\newblock {\em Math. Z. 243}, 2 (2003), 291--331.

\bibitem{maxwell}
{\sc Maxwell, J.~C.}
\newblock On the dynamical theory of gases.
\newblock {\em Philos. Trans. Roy. Soc. London Ser. A 157\/} (1867), 49–88.

\bibitem{MISCHLER:2010:HAL-00447988:1}
{\sc {M}ischler, S., and {M}ouhot, C.}
\newblock {Q}uantitative uniform in time chaos propagation for {B}oltzmann
  collision processes.
\newblock ArXiv preprint 1001.2994 (2010).

\bibitem{MR2264623}
{\sc Mischler, S., and Mouhot, C.}
\newblock Cooling process for inelastic {B}oltzmann equations for hard spheres.
  {II}. {S}elf-similar solutions and tail behavior.
\newblock {\em J. Stat. Phys. 124}, 2-4 (2006), 703--746.

\bibitem{MR2254617}
{\sc Mouhot, C.}
\newblock Explicit coercivity estimates for the linearized {B}oltzmann and
  {L}andau operators.
\newblock {\em Comm. Partial Differential Equations 31}, 7-9 (2006),
  1321--1348.

\bibitem{MR2197542}
{\sc Mouhot, C.}
\newblock Rate of convergence to equilibrium for the spatially homogeneous
  {B}oltzmann equation with hard potentials.
\newblock {\em Comm. Math. Phys. 261}, 3 (2006), 629--672.

\bibitem{MR2301289}
{\sc Mouhot, C.}
\newblock Quantitative linearized study of the {B}oltzmann collision operator
  and applications.
\newblock {\em Comm. Math. Sci.}, suppl. 1 (2007), 73--86.

\bibitem{Mouhot-Neumann}
{\sc Mouhot, C., and Neumann, L.}
\newblock Quantitative perturbative study of convergence to equilibrium for
  collisional kinetic models in the torus.
\newblock {\em Nonlinearity 19}, 4 (2006), 969--998.

\bibitem{MRS}
{\sc Mouhot, C., Russ, E., and Sire, Y.}
\newblock Fractional {P}oincar\'e inequalities for general measures.
\newblock To appear in {\it J. Math. Pures. Appl.}

\bibitem{MR2322149}
{\sc Mouhot, C., and Strain, R.~M.}
\newblock Spectral gap and coercivity estimates for linearized {B}oltzmann
  collision operators without angular cutoff.
\newblock {\em J. Math. Pures Appl. (9) 87}, 5 (2007), 515--535.

\bibitem{Pao:TSnoncutoff:74}
{\sc Pao, Y.~P.}
\newblock Boltzmann collision operator with inverse-power intermolecular
  potentials. {I}, {II}.
\newblock {\em Comm. Pure Appl. Math. 27\/} (1974), 407--428; ibid. 27 (1974),
  559--581.

\bibitem{TO}
{\sc Toscani, G.}
\newblock Entropy production and the rate of convergence to equilibrium for the
  {F}okker-{P}lanck equation.
\newblock {\em Quart. Appl. Math. 57}, 3 (1999), 521--541.

\bibitem{TV1}
{\sc Toscani, G., and Villani, C.}
\newblock Sharp entropy dissipation bounds and explicit rate of trend to
  equilibrium for the spatially homogeneous {B}oltzmann equation.
\newblock {\em Comm. Math. Phys. 203}, 3 (1999), 667--706.

\bibitem{TV2}
{\sc Toscani, G., and Villani, C.}
\newblock On the trend to equilibrium for some dissipative systems with slowly
  increasing a priori bounds.
\newblock {\em J. Stat. Phys. 98}, 5-6 (2000), 1279--1309.

\bibitem{MR1650006}
{\sc Villani, C.}
\newblock On a new class of weak solutions to the spatially homogeneous
  {B}oltzmann and {L}andau equations.
\newblock {\em Arch. Rational Mech. Anal. 143}, 3 (1998), 273--307.

\bibitem{MR1715411}
{\sc Villani, C.}
\newblock Regularity estimates via the entropy dissipation for the spatially
  homogeneous {B}oltzmann equation without cut-off.
\newblock {\em Rev. Mat. Iberoamericana 15}, 2 (1999), 335--352.

\bibitem{villani-hdr}
{\sc Villani, C.}
\newblock {\em Contribution \`a l'\'etude math\'ematique des collisions en
  th\'eorie cin\'etique (HDR)}.
\newblock PhD thesis, Univ. Paris Dauphine, France, 2000.

\bibitem{Vi:hb}
{\sc Villani, C.}
\newblock A review of mathematical topics in collisional kinetic theory.
\newblock In {\em Handbook of mathematical fluid dynamics, {V}ol. {I}}.
  North-Holland, Amsterdam, 2002, pp.~71--305.

\bibitem{MR1964379}
{\sc Villani, C.}
\newblock Cercignani's conjecture is sometimes true and always almost true.
\newblock {\em Comm. Math. Phys. 234}, 3 (2003), 455--490.

\bibitem{Vill:tspt}
{\sc Villani, C.}
\newblock {\em Topics in optimal transportation}, vol.~58 of {\em Graduate
  Studies in Mathematics}.
\newblock American Mathematical Society, Providence, RI, 2003.

\bibitem{WCUh:LBE:70}
{\sc Wang~Chang, C.~S., Uhlenbeck, G.~E., and de~Boer, J.}
\newblock In {\em Studies in Statistical Mechanics, Vol. V}. North-Holland,
  Amsterdam, 1970.

\bibitem{Wenn:stab:93}
{\sc Wennberg, B.}
\newblock Stability and exponential convergence in {$L\sp p$} for the spatially
  homogeneous {B}oltzmann equation.
\newblock {\em Nonlinear Anal. 20}, 8 (1993), 935--964.

\bibitem{MR1450762}
{\sc Wennberg, B.}
\newblock Entropy dissipation and moment production for the {B}oltzmann
  equation.
\newblock {\em J. Stat. Phys. 86}, 5-6 (1997), 1053--1066.

\end{thebibliography}

\begin{flushleft} \signld \end{flushleft}
\begin{flushleft} \signcm \end{flushleft}
\begin{flushleft} \signcv \end{flushleft}
\end{document}